\documentclass{amsart}

\usepackage{amsmath}
\usepackage{amssymb}

\DeclareMathOperator{\dom}{dom}
\newcommand{\Ord}{\textrm{Ord}}

\title{Convergence on surreals}
\author{Istv\'an Mez\H{o}}
\address{Escuela Polit\'ecnica Nacional\\Departamento de Matem\'atica\\Quito, Ecuador\\Ladr\'on de Guevara E11-253, Quito, Ecuador}
\email{istvan.mezo@epn.edu.ec}

\begin{document}

\begin{abstract}
Using the sign expansion of the surreal numbers, we give a possible notion of convergence for surreal sequences.
\end{abstract}

\maketitle

\section{The surreal numbers}

We introduce the surreal numbers via their sign expansion. A surreal number $X$ is a function on an ordinal number with range $\{+,-\}$. The class of surreals is commonly denoted by \texttt{No}.

For example, $X:\omega\to\{+,-\}$, $X(i)=+$ ($i=0,1,\dots$) denotes ''the'' infinite number $\omega=+^\omega$, while $X(0)=+$, $X(i)=-$ ($i=1,2,\dots$) denotes the infinitesimal number $\varepsilon=+-^\omega$. The natural numbers have surreal representations as $+^n$ and the dyadic rationals can be represented on the obvious way.

Some other specific surreals are
\[\omega+1=+^\omega+,\quad2\omega=+^\omega+^\omega,\quad\frac{\omega}{2}=+^\omega-^\omega,\]
and so on.

It is not obvious at all, how should we interpret limits of surreal sequences. For example, the surreal sequence $X_n=n$ should be tend to $\omega$, but, there are many numbers between the integers and $\omega$. We cite J. H. Conway, the ,,inventor'' of the surreals \cite{Conway}:

''For instance, the limit of the sequence $0,\frac{1}{2},\frac{2}{3},\frac{3}{4},\dots$ ($\omega$ terms) is not $1$, at least in the ordinary sense, because there are plenty of numbers in between. A simpler, but sometimes less convincing, example of the same phenomenon is given by the sequence $0,1,2,3,\dots$ of all finite ordinals, which one would expect tend to $\omega$, but which obviously can't since there is a whole Host of numbers greater than every finite integer but less than $\omega$. For the author's amusement, we recall some of the simplest of them:
\[\omega-1,\;\omega/2,\;\sqrt{\omega},\;\omega^{1/\omega},\dots\mbox{''}\]
That these numbers are indeed \emph{between} the sequence $0,1,2,3,\dots$ and $\omega$ can be seen by their sign expansions:
\begin{align*}
\omega-1&=+^\omega-,\\
\frac{\omega}{2}&=+^\omega-^\omega,\\
\sqrt{\omega}&=+^\omega\underbrace{-^\omega-^\omega\cdots-^\omega}_{\omega\textrm{ times!}}
\end{align*}
To put the meaning of ''between'' in a more precise form, we need the notion of ordering of surreals.

\section{Ordering and addition of the surreal numbers}

We say that $X$ is \textit{simpler} than $Y$ if $\dom X<\dom Y$ and $X(\alpha)=Y(\alpha)$ for all $\alpha<\dom X$. Here $\dom X$ is the domain of the surreal number $X$ as a function.

For example, $X=++--+-$ is simpler than $Y=++--+-++$.

Moreover, $X<Y$ if one of the following three statements hold:
\begin{enumerate}
	\item $X$ is simpler than $Y$ and $Y(\dom X)=+$,
	\item $Y$ is simpler than $X$ and $X(\dom Y)=-$,
	\item there exists a surreal number $Z$ such that $Z$ is simpler than $X$ and $Y$ and, in addition, $X(\dom Z)=-$, and $Y(\dom Z)=+$.
\end{enumerate}

So, for example $-++$ is less than $-+++$ or $-+++--+$ by the first statement; $--+-$ and $--+-++-$ are less than $--+$, since the second statement holds, and $++-+--<++++--$ by the third statement (with $Z=++$).

The notations -- like $\sqrt{\omega}$ -- are not just occasional, the surreal addition and multiplication gives sense to these expressions. We introduce addition now -- using the ordering defined above.

Let
\begin{align*}
L(X)&=\{X\restriction_{\alpha}\;:\;\alpha<\dom X\mbox{ and }X(\alpha)=+\},\\
R(X)&=\{X\restriction_{\alpha}\;:\;\alpha<\dom X\mbox{ and }X(\alpha)=-\}
\end{align*}
be the left and right set of $X$, respectively. They determine $X$ uniquely.

For instance, if $X=++-++--$, then
\begin{align*}
L(X)&=\{+,++-,++-+\},\\
R(X)&=\{++,++-++,++-++-\}.
\end{align*}
(Then one can recover $X$ concatenating a '$+$' to the greatest (last) element of $L(X)$ or concatenating a '$-$' to the greatest (last) element of $R(X)$, depending on the last sign of $X$. However, this does not work for infinite sign expansions.)

Then the sum of two surreals is defined inductively. The left and right set of $X+Y$ is \cite{Gonshor}
\begin{align*}
L(X+Y)&=\{Z+Y\;:\;Z\in L(X)\}\cup\{X+Z\;:\;Z\in L(Y)\},\\
R(X+Y)&=\{Z+Y\;:\;Z\in R(X)\}\cup\{X+Z\;:\;Z\in R(Y)\}.
\end{align*}
The multiplication defines similarly, but we shall not need it.

\section{A possible limit notion for surreal sequences}

As Conway said, it seems unnatural to pick out a special infinite surreal from all the infinite surreals saying that this is \emph{the} limit of the sequence $X_n=n$. On the other hand, taking the sign expansion into account, it is not natural to involve infinite surreals \emph{longer than} $\omega$. If we restrict us to some limit ordinals and numbers of this length, the limit becomes unique (if it exists). To be more precise, we introduce the \emph{birthday set} of a sequence $X_n$:
\[B_{X_n}:=\{\dom X_1,\dom X_2,\dom X_3,\dots\},\]
while the \emph{limit birthday} of $X_n$ is the least limit ordinal greater than or equal to $\limsup B_{X_n}$:
\[b_{X_n}:=\inf\{\alpha\in\Ord\;:\;\alpha\ge\limsup B_{X_n}\mbox{ and }\alpha\mbox{ is a limit ordinal}\}.\]
Since $B_{X_n}$ is a countably infinite set of ordinals, its limit superior exists and it will be an ordinal. Hence $b_{X_n}$ also exists.

That why we have chosen limit ordinals is obvious. On the other hand, taking the limsup of birthdays prevents us to deal with surreals ''far'' from our sequence in context. If we would take all the birthdays (lengths) of $X_n$'s, a sequence like $(1,\omega^\omega,2,3,4,5,\dots)$ would involve larger ordinals than we need.

To mention just a few example, we note that the sequence $X_n=n$ has limit birthday $\omega$, moreover $b_{(X_n=1)}=\omega$, $b_{(X_n=n\omega)}=\omega^\omega$, while $b_{(X_n=\omega/2^n)}=2\omega$.

We advice the next limit notion for surreal sequences.

\textbf{Definition.} Let $X_n$ be a sequence and let $b_{X_n}$ be its limit birthday. Moreover, let
\[H_{X_n}:=\{X\in\texttt{No}\;:\;\dom X\le b_{X_n}\}.\]
We say that $X_n$ is \emph{convergent}, if there exists an $X\in H_{X_n}$ such that for all $\alpha<b_{X_n}$ there exists an $n_0$ such that $(X_n)\restriction_{\beta_n}=X\restriction_{\beta_n}$ for all $n>n_0$. Here $\beta_n=\min\{\alpha,\dom X_n,\dom X\}$.

In this case we define $X$ as the \emph{limit} of $X_n$ and it is denoted by $\lim_{n\to\infty}X_n$.

\section{Examples}

Now let us take $X_n=n$ as a first example. The limit birthday of this sequence is $\omega$, so $H_{X_n}$ contains the surreals which have sign expansion no longer than $\omega$. If $\alpha<\omega$, say $m$, then there exists an $n_0$, namely $m$, such that $(X_n)\restriction_m=\omega\restriction_m$, if $n>m$. (Now $\beta_n=\min\{\alpha(=m),\dom X_n,\dom\omega\}=m$.)

Let us consider a bit more involved example: $X_n=1-\frac{1}{2^n}$. This sequence also has limit birthday $\omega$, so there are two possible limits: $1$ and $1-\varepsilon$, where $\varepsilon=+-^\omega$. One can see by the addition rule that $1-\varepsilon=+-+^\omega$. Since $X_n=+-+^{n-1}$, it is obvious that the limit
\[\lim_{n\to\infty}\left(1-\frac{1}{2^n}\right)=1-\varepsilon.\]

We remark that this limit notion is not additive: $X_n=Y_n=n$ tend to $\omega$, but $X_n+Y_n=2n\to\omega\neq\omega+\omega=2\omega$. The similar is true for multiplication.

Finally, we present an other example, which shows that the behaviour of this limit notion can give different results than the ordinary ''real limit''. Namely, the sequence $X_n=\frac{(-1)^n}{2^n}$ is not convergent in our context, because the sign expansions start with
\[+,-+,+--,-+++,\dots\]
The first sign alternates, hence we can not find a unique $X$ as a limit.

What is about the example of Conway? That is, what is the surreal limit of the sequence $X_n=\frac{n}{n+1}$? It is easy to see, that $X_n\to1-\varepsilon$. This is so, because for any $n>1$ there is an $m>1$ such that $1-\frac{1}{2^{m-1}}<\frac{n}{n+1}<1-\frac{1}{2^m}$. Therefore $X_n$ has the same limit as $1-\frac{1}{2^m}$.

\section{Infinite sums}

The unusual behavior of surreal limit continues to hold when we consider infinite sums. We identify an infinite sum to the limit of the sequence of its partial sums -- if the limit exists. In this case it is straightforward to see that
\[\sum_{n=1}^\infty1=\omega,\]
since tha partial sums are the natural numbers. To take another example, the harmonic series has $\infty$ as an ordinary limit, so it is convergent in the set of extended real numbers. We might think that we can find a surreal number to which $\sum_{k=1}^\infty\frac{1}{k}$ converges (a natural candidate could be $\omega$.) In contrary, this sum is divergent in the "surreal sense". This can be heuristically seen: adding $\frac{1}{n+1}$ to the partial sum $\sum_{k=1}^n\frac{1}{k}$, the integer part steadily grows, but the fractional part does not "stabilize". Hence $\sum_{k=1}^\infty\frac{1}{k}$ is divergent.

\end{document}